\documentclass[12pt,reqno]{article}

\usepackage[usenames]{color}
\usepackage{amssymb}
\usepackage{graphicx}
\usepackage{amscd}

\usepackage[colorlinks=true,
linkcolor=webgreen,
filecolor=webbrown,
citecolor=webgreen]{hyperref}

\definecolor{webgreen}{rgb}{0,.5,0}
\definecolor{webbrown}{rgb}{.6,0,0}

\usepackage{color}
\usepackage{fullpage}
\usepackage{float}

\usepackage{graphics,amsmath,amssymb}
\usepackage{amsthm}
\usepackage{amsfonts}
\usepackage{latexsym}
\usepackage{epsf}

\usepackage{verbatim}
\usepackage{enumerate}

\usepackage{tikz}
\usepackage{tcolorbox}
\usepackage[blocks]{authblk}

\usepackage{algpseudocode}

\theoremstyle{plain}

\theoremstyle{definition}

\newtheorem{example}{Example}

\theoremstyle{remark}

\begin{document}

\title{Fibonometry and Beyond}
\author{Nikhil Byrapuram}
\author{Adam Ge}
\author{Selena Ge}
\author{Sylvia Zia Lee}
\author{Rajarshi Mandal}
\author{Gordon Redwine}
\author{Soham Samanta}
\author{Daniel Wu}
\author{Danyang Xu}
\author{Ray Zhao}
\affil{PRIMES STEP}
\author{Tanya Khovanova}
\affil{MIT}
\date{}

\maketitle

\begin{abstract}
In 2013, Conway and Ryba wrote a fascinating paper called Fibonometry. The paper, as one might guess, is about the connection between Fibonacci numbers and trigonometry. We were fascinated by this paper and looked at how we could generalize it. We discovered that we weren't the first. In this paper, we describe our journey and summarize the results.
\end{abstract}

\section{Fibonometry}

Consider the sequence of Fibonacci numbers $F_n$ defined as: $F_0 = 0$, $F_1 = 1$, and for $n > 1$, we have $F_n = F_{n-1} + F_{n-2}$. We also consider its pair sequence, Lucas numbers $L_n$ defined as: $L_0 = 2$, $L_1 = 1$, and for $n > 1$, we have the same recursion $L_n = L_{n-1} + L_{n-2}$.

These numbers are involved in many famous identities. For example, we have index addition formulae
\[2F_{m+n} = F_mL_n + L_mF_n \quad \textrm{ and } \quad 2L_{m+n} = L_mL_n + 5 F_mF_n,\]
from which we can deduce the double index formulae
\[F_{2n} = F_nL_n \quad \textrm{ and } \quad 2L_{2n} = L_n^2 + 5F_n^2.\]

These formulae bear a striking similarity to trigonometric identities, starting with the angle addition formulae:
\[ \sin(\alpha +\beta )=\sin \alpha \cos \beta +\cos \alpha \sin \beta  \quad \textrm{ and } \quad \cos(\alpha +\beta )=\cos \alpha \cos \beta -\sin \alpha \sin \beta,\]
and the double-angle formulae derived from them:
\[ \sin(2\alpha )=2\sin \alpha \cos \alpha \quad \textrm{ and } \quad \cos(2\alpha )=\cos ^{2}\alpha -\sin ^{2}\alpha,\]
not to mention the angle-negation formulae
\[\sin(-x)=-\sin x \quad \textrm{ and } \quad \cos(-x)=\cos(x)\]
that are very similar to
\[F_{-n}=-(-1)^nF_n \quad \textrm{ and } \quad L_{-n}=(-1)^nL_n.\]

These similarities are based on the mathematical phenomenon called Fibonometry. The rule to convert between these parallel identities is formulated by Conway and Ryba \cite{CR} as follows.

\begin{tcolorbox}
\textbf{Fibonometry rule.} Replace an angle $\theta = p\alpha + q\beta + r\gamma + \cdots$ with a subscript $n = pa + qb + rc + \cdots$. Then replace $\sin \theta$ with $\frac{i^nF_n}{2}$ and $\cos \theta$ with $\frac{i^nL_n}{2}$. Finally, insert a factor $-5$ for any square of sines; thus, insert $(-5)^k$ for any term that contains $2k$ or $2k+1$ sines.
\end{tcolorbox}

For example, an angle $\alpha + 2\beta + 3\gamma$ is replaced with a subscript $a+2b+3c$. We start with a simple case of $\theta = \alpha$ and without sine squares.

\begin{example}
Consider the identity $\sin(-\theta) = - \sin \theta$. It becomes $\frac{i^{-n}F_{-n}}{2} = -\frac{i^nF_n}{2}$, which is equivalent to $F_{-n} = -(-1)^nF_n$.
\end{example}

Now, here's an example with squares of sines.

\begin{example}
Consider the Pythagorean identity:
\[ \sin ^{2}\theta +\cos ^{2}\theta =1.\]
It becomes $(-5)(\frac{i^{n}F_{n}}{2})^2 + (\frac{i^{n}L_{n}}{2})^2 = 1$, which simplifies to the following famous identity:
\[ L_n^{2} -5F_n^{2} = 4(-1)^n.\]
\end{example}
The Fibonometry rule emphasizes the similarity of the identities. We replace angles with subscripts, sines with Fibonacci numbers, and cosines with Lucas numbers. Then, we perform some minor adjustments.

Here's an alternate Fibonometry rule that requires only one step and doesn't have to look at sine squares.

\begin{tcolorbox}
\textbf{One-step Fibonometry rule.} Replace an angle $\theta$ with a subscript $n$ as before. Then, replace
\[ \sin\theta \quad\textrm{with}\quad \frac{\sqrt{5}}{2}i^{n-1}F_n \quad\textrm{ and }\quad \cos\theta \quad\textrm{with}\quad \frac{1}{2}i^nL_n.\]
\end{tcolorbox}

The rules are different, but how come they both work? We invite the reader to check that the rules are equivalent as long as the number of sines in every formula term has the same parity.

Did you notice that standard trigonometry formulae satisfy this parity condition? This is because the sine function is odd. Thus, if one term has an odd number of sines and another has an even number of sines, then $\theta$ and $-\theta$ can't automatically satisfy the same trigonometric formula.

We won't provide the proof for these rules; we would rather refer the readers to the 9-part series of YouTube videos \cite{Metzler} and paper \cite{Fokkink}.

The reason these rules work is the parallelism in the following formulae for Fibonacci and Lucas numbers
\[F_{n}={\frac {\varphi ^{n}-(-\varphi )^{-n}}{\sqrt {5}}} \quad \textrm{ and } \quad L_n = \varphi ^{n}+(-\varphi )^{-n}\]
as compared to sine and cosine
\[\sin \theta = \frac{e^{i\theta} - e^{-i\theta}}{2i} \quad \textrm{ and } \quad \cos \theta = \frac{e^{i\theta} + e^{-i\theta}}{2}.\]

Do you know that there are Lucas sequences in addition to Lucas numbers $L_n$ described here? We now transition to them.

\section{Lucas sequences}

There are other families of sequences that are similar to Fibonacci and Lucas numbers. Confusingly, they are called Lucas sequences.

\subsection{Definition}

The Lucas sequences $U_{n}(P,Q)$ and $V_{n}(P,Q)$ are integer sequences that satisfy the recurrence relation $x_{n}=P x_{n-1} - Q x_{n-2}$, where $P$ and $Q$ are fixed integers. We also have initial conditions: $U_{0}(P,Q) = 0$, $U_{1}(P,Q) = 1$, $V_{0}(P,Q) = 2$, and $V_{1}(P,Q) = P$. Given two integer parameters $P$ and $Q$, sequence $U_{n}(P,Q)$ is called the Lucas sequence of the first kind, and sequence $V_{n}(P,Q)$ is called the Lucas sequence of the second kind.

The Lucas sequences are a generalization of Fibonacci and Lucas numbers. Indeed, $U_{n}(1,-1)$ are Fibonacci numbers, and $V_{n}(1,-1)$ are Lucas numbers.

The name `Lucas sequences' is very confusing, given that a sequence of Lucas numbers exists, which becomes a Lucas sequence in two different senses. To resolve the ambiguity, we suggest calling Lucas sequences, \textit{general Lucas sequences} and Lucas numbers \textit{proper Lucas numbers}.

\begin{example}
The pairs of general Lucas sequences for fixed values of $P$ and $Q$ contain other famous pairs of sequences; for example, we get:
\begin{itemize}
\item $U_{n}(2,-1)$ are Pell numbers, and $V_{n}(2,-1)$ are Pell-Lucas numbers,
\item $U_{n}(1,-2)$ are Jacobsthal numbers, and $V_{n}(1,-2)$ are Jacobsthal-Lucas numbers,
\item $U_{n}(3,2)$ are numbers $2^n-1$, and $V_{n}(3,2)$ are numbers $2n+1$.
\end{itemize}
\end{example}

The characteristic equation of the recurrence relation for general Lucas sequences $U_{n}(P,Q)$ and $V_{n}(P,Q)$ is:
\[x^{2}-Px+Q=0.\]
Its discriminant $D=P^{2}-4Q$ plays an important role in this paper and elsewhere.

\subsection{Generalized Fibonometry rule}

Can the Fibonometry rule be generalized to general Lucas sequences? First, we need an additional notation. Let $\ell=\frac{P+\sqrt{D}}{2}$. Now, we can present the rule \cite{Metzler}.

\begin{tcolorbox}
\textbf{Beyond Fibonometry rule.} Replace an angle $\theta$ with a subscript $n$ as before. Then, replace
\[\sin\theta \textrm{ with } \frac{\sqrt{D}}{2i}\left(\frac{-1}{\sqrt{Q}}\right)^nU_{n} \quad\textrm{and}\quad \cos\theta \textrm{ with } \frac{1}{2}\left(\frac{-1}{\sqrt{Q}}\right)^nV_{n}.\]
\end{tcolorbox}

Now, let's look at some examples of the application of this rule.
\begin{example}
Consider the identity $\sin\alpha+\sin\beta=2\sin\left(\frac{\alpha+\beta}{2}\right)\cos\left(\frac{\alpha-\beta}{2}\right)$. Applying the rule, we get
\[\frac{\sqrt{D}}{2i}\left(\frac{-1}{\sqrt{Q}}\right)^mU_{m} + \frac{\sqrt{D}}{2i}\left(\frac{-1}{\sqrt{Q}}\right)^nU_{n} = 2 \frac{\sqrt{D}}{2i}\left(\frac{-1}{\sqrt{Q}}\right)^{\frac{m+n}{2}}U_{\frac{m+n}{2}} \frac{1}{2}\left(\frac{-1}{\sqrt{Q}}\right)^{\frac{m-n}{2}}V_{\frac{m-n}{2}},\]
which can be simplified to
\[(-1)^m(\sqrt{Q})^nU_m + (-1)^n(\sqrt{Q})^mU_n = (-1)^m(\sqrt{Q})^nU_{\frac{m+n}{2}}V_{\frac{m-n}{2}}.\]

Dividing each term by $(-1)^m(\sqrt{Q})^n$, we get
\[ U_m+(-\sqrt{Q})^{m-n}U_n = U_{\frac{m+n}{2}}V_{\frac{m-n}{2}}. \]
Since $m$ and $n$ must be of the same parity in order for $\frac{m+n}{2}$ to be an integer, we can assume that $m$ and $n$ have the same parity and simplify the above equation further
\[ U_m+Q^{\frac{m-n}{2}}U_n=U_{\frac{m+n}{2}}V_{\frac{m-n}{2}}. \]
\end{example}

\begin{example}
The trigonometry identity
\[ \sin \alpha \sin(\beta-\gamma) + \sin \beta \sin(\gamma-\alpha) + \sin \gamma \sin(\alpha-\beta) = 0, \]
is equivalent to
\[Q^nU_\ell U_{m-n} + Q^{\ell}U_mU_{n-\ell} + Q^{m}U_nU_{\ell-m} = 0.\]
\end{example}

Table~\ref{tab:FLGeneral} shows examples of conversions. Keep in mind that in the bottom two rows, we assume that $m-n$ is even.

\begin{table}[ht!]
\begin{center}
\begin{tabular}{|c|c|}
\hline
Trigonometry					& General Lucas Sequences                \\ \hline
$\sin^2 \alpha + \cos^2 \alpha = 1$		& $V_n^2-DU_n^2=4Q^n$                  \\ \hline
$\sin(-\alpha) = - \sin \alpha$		& $-Q^nU_{-n} =U_n$               \\ \hline
$\cos(-\alpha)=\cos\alpha$			&$V_n=Q^nV_{-n}$\\\hline
$\sin 2\alpha = 2\cos \alpha \sin \alpha$	& $U_{2n}=U_nV_n$                     \\ \hline
$\cos 2\alpha = 2\cos^2 \alpha - 1$		& $V_{2n}=V_n^2-2Q^n$                 \\ \hline
$\cos 2\alpha=1-2\sin^2\alpha$		&$V_{2n}=2Q^n-DU^2_n$\\\hline
$\sin 3\alpha = 3\sin\alpha - 4\sin^3\alpha$	& $U_{3n} = 3Q^nU_n + DU^3_n$	\\ \hline
$\cos 3\alpha=4\cos^3\alpha-3\cos\alpha$	&$V_{3n}=V^3_n-3Q^nV_n$\\\hline
$\sin(\alpha+\beta)=\sin\alpha\cos\beta + \cos\alpha\sin\beta$	& $2U_{m+n}=U_mV_n+U_nV_m$ 	\\ \hline
$\cos(\alpha+\beta)=\cos\alpha\cos\beta - \sin\alpha\sin\beta$	& $2V_{m+n}=V_mV_n+ DU_mU_n$ \\ \hline
$\sin\alpha\sin\beta=\frac{\cos(\alpha-\beta)-\cos(\alpha+\beta)}{2}$	& $DU_mU_n = V_{m+n}- Q^nV_{m-n}$	\\ \hline
$\cos\alpha\cos\beta=\frac{\cos(\alpha+\beta)+\cos(\alpha-\beta)}{2}$	&$V_mV_n=V_{m+n}+Q^nV_{m-n}$\\\hline
$\sin\alpha\cos\beta=\frac{\sin(\alpha+\beta)+\sin(\alpha-\beta)}{2}$	&$U_mV_n=U_{m+n}+Q^nU_{m-n}$\\\hline
$\sin\alpha+\sin\beta=2\sin\left(\frac{\alpha+\beta}{2}\right)\cos\left(\frac{\alpha-\beta}{2}\right)$	& $U_m+Q^{\frac{m-n}{2}}U_n=U_{\frac{m+n}{2}}V_{\frac{m-n}{2}}$	\\ \hline
$\cos\alpha-\cos\beta=-2\sin(\frac{\alpha+\beta}{2})\sin(\frac{\alpha-\beta}{2})$&$V_m-Q^{\frac{m-n}{2}}V_n=DU_{\frac{m+n}{2}}U_{\frac{m-n}{2}}$\\\hline
\end{tabular}
\end{center}
\caption{Trigonometry identities converted to general Lucas sequences}
\label{tab:FLGeneral}
\end{table}

\begin{example}
Let us consider $U_n(1,-2)$ (Jacobsthal numbers) and $V_n(1,-2)$ (Jacobsthal-Lucas numbers) corresponding to $P=1$, $Q=-2$. The first few terms of these sequences are shown as rows 2 and 3 in Table~\ref{tab:Jacobsthal}. The equation $\sin^2 \alpha + \cos^2 \alpha = 1$ corresponds to $V^2_n-9U^2_n = 4Q^n$. The value for the expression $V^2_n-9U^2_n$ is shown as row 4 in Table~\ref{tab:Jacobsthal}. The reader can verify that it is equal to $4Q^n$.

\begin{table}[ht!]
\centering
\begin{tabular}{|c|c|c|c|c|c|c|c|}
\hline
$n$ 		& 0 & 1 	& 2 & 3      & 4 & 5 & 6 \\ \hline
$U_n(1,-2)$ 	& 0 & 1 	& 1 & 3      & 5 & 11 & 21 \\ \hline
$V_n(1,-2)$ 	& 2 & 1 	& 5 & 7      & 17 & 31 & 65 \\ \hline
$V^2_n-9U^2_n$ & 4 & $-8$ 	& 16 & $-32$ & 64 & $-128$ & 256 \\ \hline
\end{tabular}
\caption{The first few $U_n(1,-2)$ and $V_n(1,-2)$ numbers.}
\label{tab:Jacobsthal}
\end{table}
\end{example}

In summary, generalized Fibonometry is a powerful tool---it allows us to remember only one set of identities to get the other set.

\subsection{A tricky example}

Given $w, x, y, z$ such that $w+z=x+y$, the function $f(t)$ defined as
\[\sin(w + t)\sin(z + t) - \sin(x + t)\sin(y + t)\]
is a constant independent of $t$. The reader can verify it using derivatives or converting $f(t)$ into
\[\frac{1}{2}(\cos(w-z)-\cos(w+z+2t))- \frac{1}{2}(\cos(x-y)+ \cos(x+y+2t)).\]
The cosines of $w+z+2t$ and $x+y+2t$ cancel out because $w+z+2t=x+y+2t$. Now we have $f(t)=\frac{1}{2}(\cos(w-z)-\cos(x-y))$, which is an expression without $t$, so it is a constant function.

By using Fibonometry, this example converts to the following Fibonacci identity. Given integers $a, b, c, d$ such that $a+d=b+c$, the function $g(r)$ defined as
\[(-1)^r(F_{a+r}F_{d+r} - F_{b+r}F_{c+r})\]
is a constant independent of $r$.

Other identities can be derived from this. For example, plugging in $a=d=n$, $b=n+m$, and $c=n-m$ converts the function above into another constant function independent of $r$:
\[(-1)^r(F_{n+r}F_{n+r} - F_{n+m+r}F_{n-m+r}).\]
Thus, we can equate the values for $r = 0$ and $r = m-n$, getting
\[(F_{n}F_{n} - F_{n+m}F_{n-m}) = (-1)^{m-n}(F_{m}F_{m} - F_{2m}F_{0}),\]
which is equivalent to Catalan's identity: 
\[F^2_n - F_{n-m}F_{n+m} = (-1)^{m-n}F^2_m,\]
which, in turn, can be converted into Cassini's identity with $m=1$, yielding $F_{n-1}F_{n+1}-F_{n}^{2}=(-1)^{n}$.

The example above was introduced in \cite{Metzler}. We can expand it to Lucas sequences as follows.
Similarly, we can derive the counterpart for the general Lucas sequence as follows. Given integers $a, b, c, d$ such that $a+d=b+c$, the function $g(r)$ defined as
\[Q^{-r}(U_{a+r}U_{d+r} - U_{b+r}U_{c+r})\]
is a constant independent of $r$. We then can derive the counterpart of the Catalan's identity:
\[ U^2_n-U_{n-m}U_{n+m}=Q^{n-m}U^2_m.\]

\section{Hyperbolic Case}

Unsurprisingly, trigonometric identities are similar to trigonometric identities for hyperbolic functions.

\subsection{Hyperbolic trigonometric functions}

Hyperbolic trigonometric functions $\cosh$ and $\sinh$ are defined using the following formulae:
\[\cosh \theta = \frac{e^{\theta} + e^{-\theta}}{2} \quad \textrm{ and } \quad \sinh \theta = \frac{e^{\theta} - e^{-\theta}}{2i}.\]

The connection to trigonometric functions can be described as follows
\[\cosh \theta = \cos{\frac{\theta}{i}} \quad \textrm{ and } \quad \sinh \theta = i\sin{\frac{\theta}{i}}.\]

It isn't surprising that there exists a conversion rule \cite{Osborn}.

\begin{tcolorbox}
\textbf{Osborn's Rule}. A trigonometric identity can be converted to an analogous identity for hyperbolic functions by exchanging trigonometric functions with their hyperbolic counterparts and then flipping the sign of each term involving the product of two hyperbolic sines.
\end{tcolorbox}

\begin{example} The expression
\[\cos(\alpha-\beta) = \cos \alpha \cos \beta + \sin \alpha \sin \beta\] converts to 
\[\cosh(\alpha-\beta) = \cosh \alpha \cosh \beta - \sinh \alpha \sinh \beta.\]
\end{example}

This rule is parallel to the Fibonometry rule. Specifically, it has a separate clause for squares. We can have a one-step rule as follows. Replace $\sin \theta$ with $i\sinh \theta$ and $\cos \theta$ with $\cosh \theta$.

\subsection{Converting hyperbolic trigonometric identities to Fibonacci}

Can we convert hyperbolic trigonometric identities into Fibonacci-Lucas identities? Yes, we can. This rule appeared in Vajda's book \cite{V} in 2008, way before the Conway-Ryba paper.

\begin{tcolorbox}
\textbf{HyperFibonometry rule.} Replace an angle $\theta$ with a subscript $n$ as before. Then, replace
\[\sinh \theta \quad \textrm{with} \quad \frac{F_{n}\sqrt{5}}{2i^{n}} \quad\textrm{and}\quad \cos\theta \quad \textrm{with} \quad \frac{L_{n}}{2i^{n}}.\]
\end{tcolorbox}

This rule doesn't mention the squares of hyperbolic sines. The appeal of the rule with squares is that it is more visual, so we add such a rule: For any hyperbolic function identity, replace $\cosh \theta$ with $\frac{i^n L_n}{2}$ and $\sinh$ with $\frac{i^n F_n}{2}$. For any square of hyperbolic sines, insert a factor of $5$.

The rule is exactly the same as the Fibonometry rule, except that we use 5 instead of $-5$.

\begin{example}
Take the hyperbolic identity
\[\cosh^2 \theta - \sinh^2 \theta = 1.\]
Applying our rule gives $(\frac{i^{n} L_n}{2})^2 - 5(\frac{i^{n} F_n}{2})^2 = 1$, which is the famous identity
\[L_n^2 - 5 F_n^2 = 4(-1)^n.\]
\end{example}

\subsection{Hyperbolic trigonometric identities and general Lucas sequences}

We can combine the rules to get a direct conversion rule from hyperbolic trigonometric identities to general Lucas sequences.

\begin{tcolorbox}
\textbf{Beyond HyperFibonometry rule.} Replace an angle $\theta$ with a subscript $n$ as before. Then, replace
\[\sinh \theta \quad \textrm{with} \quad \frac{\sqrt{D}}{2}\left(\frac{-1}{\sqrt{Q}}\right)^nU_n \quad \textrm{ and }\quad \cos\theta \quad \textrm{with} \quad \frac{1}{2}\left(\frac{-1}{\sqrt{Q}}\right)^nV_n.\]
\end{tcolorbox}

\begin{example}
We can convert $\cosh^2 \theta - \sinh^2 \theta =1$ to $V^2_n-DU^2_n=4Q^n$.
\end{example}

\section{Going Further}

Many trigonometric identities contain $\tan$ and other functions. Below, we will only consider $\tan$; however, other functions can be treated similarly. Based on the Beyond Fibonometry rule, we can have a direct conversion rule for $\tan$. We should replace
\[ \tan \theta \quad \textrm{ with } \quad \frac{\sqrt{D}}{i} \cdot \frac{U_n}{V_n}.\]

\begin{example}
We can convert the following trigonometric identity
\[ \tan \frac{\theta}{2} = \frac{1-\cos\theta}{\sin\theta} \]
into
\[\frac{\sqrt{D}}{i} \cdot \frac{U_{\frac{n}{2}}}{V_{\frac{n}{2}}} = \frac{1-\frac{1}{2}\left( \frac{-1}{\sqrt{Q}} \right)^n V_n}{\frac{\sqrt{D}}{2i}\left( \frac{-1}{\sqrt{Q}} \right)^n U_n},\]
which eventually gives the identity
\[ DU_{\frac{n}{2}}U_n = V_{\frac{n}{2}}V_n-2Q^{\frac{n}{2}}V_{\frac{n}{2}}\]
for even $n$.
\end{example}

\begin{example}
The trigonometric identity
\[\tan(x+y) = \frac{\tan x+\tan y}{1-\tan x \tan y}.\]
It can be converted to
\[\frac{U_{m+n}}{V_{m+n}}=\frac{U_mV_n+U_nV_m}{V_mV_n+DU_mU_n}.\]
Similarly, the trigonometric identity
\[\cos(2x)=\frac{1-\tan^2(x)}{1+\tan^2(x)}\]
can be converted to
\[V_{2n}=2Q^n\cdot \frac{V^2_n+DU^2_n}{V^2_n-DU^2_n}.\]
\end{example}

\section{Going Backwards}

Of course, we can go backward, starting with a Fibonacci identity and reverting it to a trigonometric identity.

\begin{tcolorbox}
\textbf{Backwards rule.} Replace a subscript $n$ with an angle $\theta$ similar to before. Then, replace
\[F_n \textrm{ with } \frac{2\sin(\theta)}{\sqrt{5}i^{n-1}} \quad \text{ and } \quad L_n \textrm{ with } \frac{2\cos(\theta)}{i^{n}}.\]
\end{tcolorbox}

\begin{example}
Consider the following identity from \cite{Keskin}
\[(-1)^{n+m}L^2_{n+m}-5(-1)^nF^2_n-5(-1)^mF^2_m = 5(-1)^{n+m}F_nF_mL_{n+m}+4.\]
Using the rule above, it becomes:
\[\cos^2(\alpha+\beta)+\sin^2\alpha+\sin^2\beta+2\sin \alpha\sin\beta\cos(\alpha+\beta)=1.\]
\end{example}

\section{Conclusion}

Everything in the world is connected, especially in mathematics. We started with Fibonometry, a branch of mathematics that connects two seemingly unrelated topics: trigonometry with Fibonacci and Lucas numbers. We extended it further to hyperbolic trigonometry on one hand and to general Lucas sequences on the other hand.

So, the next time you're pondering the beauty of a sine wave, remember it might just be whispering the secrets of the Fibonacci sequence in your ear. Or it might hold patterns waiting to be discovered.

\section{Acknowledgments}

We are grateful to the MIT PRIMES STEP program for allowing us the opportunity to do this research.


\begin{thebibliography}{9}

\bibitem{CR} John Conway and Alex Ryba, Fibonometry. \textit{The Mathematical Gazette} 97 (2013), pp.~494--495.

\bibitem{Fokkink} Robert Fokkink, The Pell Tower and Ostronometry. arXiv preprint arXiv:2309.01644 (2023).

\bibitem{Keskin} Refik Keskin, Three identities concerning Fibonacci and Lucas numbers, \textit{Notes on Number Theory and Discrete Mathematics}, 20:5, (2014), pp.~44--48.

\bibitem{Metzler} David Metzler, Fibonacci Numbers and Complex Trigonometry. YouTube series, (2015) \url{https://www.youtube.com/playlist?list=PLP8japRypwHmtdpyApWY_R6A_c0BKbT1W}.

\bibitem{Osborn} G.\ Osborn, 109.[D. 6. d.] Mnemonic for hyperbolic formulae. \textit{The Mathematical Gazette} 2(34), (1902) p.~189.

\bibitem{V} Steven Vajda, \textit{Fibonacci and Lucas numbers, and the golden section: theory and applications}. Courier Corporation, 2008.

\end{thebibliography}
\end{document}